\documentclass{amsart}
\usepackage{amsmath}
\usepackage{amssymb}
\usepackage{amsthm}
\usepackage[all]{xy}
\usepackage[mathscr]{eucal}
\usepackage{color}

\newtheorem{theorem}{Theorem}[section]
\newtheorem{proposition}[theorem]{Proposition}

\newtheorem{lemma}[theorem]{Lemma}
\theoremstyle{definition}
\newtheorem{definition}[theorem]{Definition}

\theoremstyle{remark}
\newtheorem{remark}[theorem]{\bf Remark}

\makeatother

\newcommand{\onto}[1]{\stackrel{#1}{\to}}
\newcommand{\Gal}{\mathrm{Gal}}
\newcommand{\inj}{\hookrightarrow}
\newcommand{\surj}{\twoheadrightarrow}
\newcommand{\cO}{\mathcal{O}}
\newcommand{\fa}{\mathfrak{a}}
\newcommand{\sn}{\smallskip\noindent}
\newcommand{\bn}{\bigskip\noindent}

\newcommand{\Pm}{(\mathrm{P}_m)}
\newcommand{\PPm}{(\mathrm{P}'_m)}
\newcommand{\QQm}{(\mathrm{Q}'_m)}
\newcommand{\PPPm}{(\mathrm{P}''_m)}
\newcommand{\Qm}{(\mathrm{Q}''_m)}
\newcommand{\ii}{\textbf{i}}
\newcommand{\uu}{\textbf{u}}
\newcommand{\ff}{\widetilde{\varphi}}
\newcommand{\fp}{\widetilde{\psi}}
\newcommand{\fe}{\mathcal{FE}}
\newcommand{\bs}{\mathcal{S}}
\newcommand{\ssb}{\mathscr{B}}

\title[Ramification of local fields and Fontaine's property $\Pm$]
{Ramification of local fields and\\ Fontaine's property $\Pm$}
\author{Manabu Yoshida}
\address{
	Graduate School of Mathematics, Kyushu University,
	Fukuoka 819-0395, Japan
}
\email{m-yoshida@math.kyushu-u.ac.jp}

\thanks{The author is supported in part by JSPS Core-to-Core Program 18005.}

\begin{document}
\maketitle
\begin{abstract}
We prove that 
the ramification filtration of the absolute Galois group of 
a complete discrete valuation field with perfect residue field 
is characterized in terms of Fontaine's property $\Pm$. 
\end{abstract}
\section{Introduction}
\quad Let $K$ be a complete discrete valuation field with perfect residue 
field $k$ 
of characteristic $p>0$, 
$\cO_K$ its valuation ring, $v_K$ its valuation normalized by $v_K(K^{\times})=\mathbb{Z}$, 
$K^{\mathrm{alg}}$ a fixed algebraic closure of $K$ and $\bar{K}$ the separable 
closure of $K$ in $K^{\mathrm{alg}}$. 
In this paper, we construct a certain decreasing 
filtration of the absolute Galois group $G_K:=\Gal(\bar{K}/K)$ 
to measure 
the ramification of extensions of $K$. 
If $E$ is an algebraic extension of $K$, we denote by $\cO_E$ the integral 
closure of $\cO_K$ in $E$. 
The valuation $v_K$ can be extended to $E$ uniquely and the extension is 
also denoted by $v_K$. 
For an algebraic extension $E$ of $K$ and a real number $m$, 
we put $\fa_{E/K}^m:=\{x \in \cO_{E}\ |\ v_K(x) \geq m\}$, which is an ideal of $\cO_E$.  
For a finite Galois extension $L/K$ and a real number $m$, 
we consider the following property studied in \cite{Fon85}:
\begin{quote}
\begin{itemize}
\item[$\Pm$] 
\emph{For any algebraic extension} $E/K$, \emph{if there exists an} 
$\cO_K$-\emph{algebra homomorphism} $\cO_L \to \cO_E/\fa_{E/K}^m$,
\emph{then there exists a} $K$-\emph{embedding} $L \inj E$. 
\end{itemize}
\end{quote}
For a finite Galois extension $L$ of $K$, 
we put 
\begin{center}
$m_{L/K}:=\inf \{ m \in \mathbb{R}\ |\ \Pm$ is true for 
$L/K \}$.
\end{center}
The property $\Pm$ is stable under composition of 
extensions of $K$ (Prop.\ \ref{composite}). 
Hence we can define two filtrations of $G_K$ as follows:
For a real number $m$, 
we denote by $\bar{K}_{<m}$ 
(resp.\ $\bar{K}_{\leqslant m}$) the composite field 
of all finite Galois extensions $L$ of $K$ in $\bar{K}$ such that 
$m_{L/K} < m$ (resp.\ $m_{L/K} \leq m$). 
We define two closed normal subgroups $G_{K}^{\geqslant m}$ 
and $G_{K}^{>m}$  
of $G_K$ by
\[G_{K}^{\geqslant m}:=\Gal(\bar{K}/\bar{K}_{<m}),\quad 
G_{K}^{>m}:=\Gal(\bar{K}/\bar{K}_{\leqslant m}).\]
The filtration $(G_K^{\geqslant m})_{m \in \mathbb{R}}$ satisfies 
$\bigcap_m G_{K}^{\geqslant m}=1$ and $G_{K}^{\geqslant 0}=G_K$ 
(Thm.\ \ref{properties} (i)). 
Moreover, $G_{K}^{\geqslant 1}$ is the inertia 
subgroup of $G_K$ and $G_{K}^{>1}$ 
is the wild inertia subgroup of $G_K$ (Thm.\ \ref{properties} (iii), 
Rem.\ \ref{WildInertia}). 

On the other hand, we denote by $G_K^{(m)}$ the $m$th upper numbering ramification group 
in the sense of \cite{Fon85}. 
Namely,\ we put $G_K^{(m)}:=G_K^{m-1}$, where the latter is the $m$th 
upper numbering ramification 
group defined in \cite{Serre}. 
In addition, we put $G_K^{(m+)}:= \overline{\bigcup_{m'>m} G_K^{(m')}}$, 
where the overline means the closure with respect to the Krull topology. 
These define two decreasing filtrations of $G_K$ and 
they are well-known in the classical ramification theory. 

We denote by $\bar{K}_{(m)}$ (resp.\ $\bar{K}_{(m+)}$) the fixed field 
of $\bar{K}$ by $G_K^{(m)}$ (resp.\ $G_K^{(m+)}$). 
Our main result in this paper is:
\begin{theorem}\label{Main}
For a real number $m$,
we have $\bar{K}_{<m}=\bar{K}_{(m)}$ and 
$\bar{K}_{\leqslant m}=\bar{K}_{(m+)}$,
so that 
\[G_{K}^{\geqslant m}=G_K^{(m)},\quad G_{K}^{>m}=G_K^{(m+)}.\]
\end{theorem}
We prove this theorem by showing the equality $m_{L/K}=u_{L/K}$ for a finite Galois 
extension $L$ of $K$, 
where $u_{L/K}$ is the greatest upper ramification break of $L/K$ in the sense of 
\cite{Fon85}. 

The property $\Pm$ is useful for obtaining ramification 
bounds for certain Galois 
representations 
(\cite{Caruso}, \cite{Fon85}, \cite{Hattori s}).
Indeed, Fontaine proved the following:
in the case where the characteristic of $K$ is $0$, for an integer $n \geq 1$, 
if we denote by $\mathscr{G}$ a finite flat group scheme over $\cO_K$ killed by $p^n$, then 
the ramification of $\mathscr{G}(\bar{K})$ 
is bounded by $m$ 
(meaning that $G_K^{(m)}$ acts trivially on $\mathscr{G}(\bar{K})$) 
if $m > e(n+ 1/(p-1))$, 
where $e$ is the absolute ramification index of $K$ (\cite{Fon85}, Thm.\ A). 
He obtained the ramification bound by showing that if 
$\Pm$ is true for a finite Galois extension $L/K$ and a real 
number $m$ then $m>u_{L/K}-e_{L/K}^{-1}$, 
where $e_{L/K}$ is the ramification index of $L/K$ (Prop.\ \ref{Fontaine} (ii)). 
The equality $m_{L/K}=u_{L/K}$ is a refinement of this result. 
Hattori (\cite{Hattori s}) generalized this kind of ramification bound 
to the case of semi-stable torsion representations.
Our equality was used in \cite{Hattori s}, 
Proposition 5.6, to improve his bound.


In Section 2, 
we study some properties of $\Pm$ and the number $m_{L/K}$. 
By using these results, we define our filtrations of $G_K$ 
and deduce its properties. 
In Section 3, 
after recalling the classical ramification theory for Galois extensions of $K$ 
(\cite{Fon85}, \cite{Serre}), 
we show the equality $m_{L/K}=u_{L/K}$ 
to prove Theorem $\ref{Main}$. 
In Section 4, we begin with a review of the ramification theory of Abbes 
and Saito (\cite{Abbes-Saito 1}, \cite{Abbes-Saito 2}). 
Their theory does not require the assumption that the residue field $k$ 
is perfect. 
Then we consider the property $\Pm$ in the imperfect 
residue field case, and translate our results in Section 
$\ref{RamificationBreaks}$ 
into the language of their theory. 
In the Appendix, we prove a Galois theoretic property on filtrations 
of the absolute Galois group of an arbitrary field. 
Theorem \ref{Main} is proved by the equality $m_{L/K}=u_{L/K}$ and the 
property checked in the Appendix.

\bn
$Convention$ $and$ $Notation$. 
Fix an algebraic closure $K^{\mathrm{alg}}$ of $K$ 
and denote by $\bar{K}$ the separable closure of $K$ in $K^{\mathrm{alg}}$. 
We assume throughout that all algebraic extensions of $K$ under discussion 
are contained in $K^{\mathrm{alg}}$. 
If $E$ is an algebraic extension of $K$, 
then we denote by $e_{E/K}$ the ramification index of $E/K$ and by $\cO_E$ 
the integral closure of $\cO_K$ in $E$. 
The valuation $v_K$ of $K$ extends to $K^{\mathrm{alg}}$ 
uniquely and the extension is also denote by $v_K$.

\bn
$Acknowledgments$. The author would like to express his deepest gratitude 
to his adviser 
Yuichiro Taguchi 
for introducing him to the problem, and for reading 
preliminary manuscripts of this paper carefully.  
He thanks Toshiro Hiranouchi for communicating Lemma 
$\ref{KrasnerImperfect}$ to him. 
He also thanks Seidai Yasuda for organizing the conference 
\emph{Ramification Theory in Arithmetic Geometry}, Kobe, 2009. 
It provided the author a great opportunity to receive many 
helpful comments by Shin Hattori and Yoichi Mieda 
for the proofs in Section $\ref{RamificationTheoryViaPm}$. 
He thanks Shinya Harada and Yoshiyasu Ozeki who taught 
him much about local fields. 
Finally, he is grateful to the referee for pointing out 
a mistake in Proposition \ref{pro2} and for suggesting many improvements of 
this paper. 

\section{Ramification theory via $\Pm$}\label{RamificationTheoryViaPm}
\quad In this section, we study the property $\Pm$. 
For a finite Galois extension $L$ of $K$, 
we put
\begin{center}
$m_{L/K}:=\inf \{ m \in \mathbb{R}\ |\ \Pm$ is true for 
$L/K \}$.
\end{center} 
If $L=K$, the property $\Pm$ holds for all real numbers $m$, 
so that we have $m_{L/K}=-\infty$. 
The following proposition is a basic property of the number $m_{L/K}$:
\begin{proposition}\label{pro1}
Let $L$ be a finite Galois extension of $K$ such that $L \not= K$. 
Then $\Pm$ is not true for $L/K$ and any real number $m \leq 0$, 
and is true for sufficiently large real number $m$. 
In particular, the number $m_{L/K}$ is non-negative and finite.  
\end{proposition}
\begin{proof}
For any real number $m \leq 0$, $\cO_K/\fa_{K/K}^m$ is zero ring. 
Then the zero map $\cO_L \to \cO_K/\fa_{K/K}^m$ is an $\cO_K$-algebra homomorphism. 
However, there is no $K$-embedding $L \inj K$ by assumption. 
Hence $\Pm$ is not true for $L/K$ and any real number $m \leq 0$. 
Thus we have $m_{L/K} \geq 0$. 
Next, we show that $\Pm$ is true for sufficiently large real number $m$. 
Choose an element $\alpha$ of $\cO_L$ such that $\cO_L=\cO_K[\alpha]$. 
Let $P$ be the minimal polynomial of $\alpha$ over $K$ and 
$\alpha=\alpha_1,\dots,\alpha_n$ the zeros of $P$ in $\bar{K}$. 
Suppose there exists an $\cO_K$-algebra homomorphism 
$\eta:\cO_L \to \cO_E/\fa_{E/K}^m$ for an algebraic extension $E$ of $K$ 
and $m>n \sup_{i \ne 1}v_K(\alpha-\alpha_i)$. 
Then we have $v_K(P(\beta)) \geq m$, 
where $\beta$ is a lift of $\eta(\alpha)$ in $\cO_E$. 
By the inequalities 
\[n\ \underset{i}{\sup}\ v_K(\beta-\alpha_i) \geq 
v_K(P(\beta)) \geq m > n \sup_{i \ne 1}\ v_K(\alpha-\alpha_i),\]
we have 
$v_K(\beta-\alpha_{i_0})>\sup_{i \ne 1}v_K(\alpha-\alpha_i)$ for some 
$i_0$. 
By Krasner's lemma, we have $K(\alpha_{i_0}) \subset K(\beta)$. 
Thus we obtain a $K$-embedding 
$L=K(\alpha)$ $\onto{\cong}$ $K(\alpha_{i_0})$ 
$\subset$ $K(\beta)$ $\subset$ $E$. 
Hence $\Pm$ is true for $m>n \sup_{i \ne 1}v_K(\alpha-\alpha_i)$. 
Therefore, we have $m_{L/K} \leq n \sup_{i \ne 1} v_K(\alpha-\alpha_i) 
<\infty$. 
\end{proof}

The following proposition often allows us to assume $L/K$ is totally ramified: 
\begin{proposition}\label{pro2}
Let $L$ be a finite Galois extension of $K$ 
and $K'$ an arbitrary finite separable 
extension of $K$. 
Put $e':=e_{K'/K}$. 
If $\Pm$ is true for $L/K$, 
then $\mathrm{(P}_{e'm}\mathrm{)}$ is true for $LK'/K'$. 
Moreover, if $K'/K$ is an unramified subextension of $L/K$ such that $L \not= K'$, 
then the converse is true. 
In particular, we have $m_{LK'/K'} \leq e'm_{L/K}$ with 
equality if $K'/K$ is an unramified subextension of $L/K$ such that $L \not= K'$. 
\end{proposition}
\begin{proof}
Put $L':=LK'$. 
First, we assume that 
$\Pm$ is true for $L/K$ and a real number $m$. 
Then we want to show that $(\mathrm{P}_{e'm})$ is also true for $L'/K'$. 
Suppose there exists an $\cO_{K'}$-algebra homomorphism 
$\eta:\cO_{L'} \to \cO_E/\fa_{E/K'}^{e'm}$ for an algebraic extension 
$E$ of $K'$. 
Then the composite map defined by 
\[\eta':\cO_L \inj \cO_{L'} 
\onto{\eta} \cO_E/\fa_{E/K'}^{e'm}=\cO_E/\fa_{E/K}^{m} \]
is an $\cO_K$-algebra homomorphism. 
Since $\Pm$ is true for $L/K$, 
there exists a $K$-embedding $L \inj E$ corresponding to $\eta'$. 
Since $L/K$ is a Galois extension, 
there exists a $K'$-embedding $L'=LK' \inj E$. 
Hence (P$_{e'm}$) is true for $L'/K'$. 
Next, we assume $K'/K$ is an unramified subextension of $L/K$ such that $L \not=K'$ 
and $\Pm$ is true for $L/K'$ and $m$. 
Note that $m>0$ by Proposition \ref{pro1}. 
Then we want to show that $\Pm$ is also true for $L/K$ and $m$. 
Suppose there exists an $\cO_K$-algebra homomorphism 
$\eta:\cO_L \to \cO_E/\fa_{E/K}^m$ 
for an algebraic extension $E$ of $K$. 
The composite map 
\[ \eta':\cO_{K'} \inj \cO_L \onto{\eta} \cO_E/\fa_{E/K}^m \]
is an $\cO_K$-algebra homomorphism. 
Since $K'/K$ is unramified, 
$\eta'$ lifts to an $\cO_K$-algebra homomorphism $\cO_{K'} \to \cO_E$. 
Hence $\eta$ is an $\cO_{K'}$-algebra homomorphism. 
By the property $\Pm$, 
there exists a $K'$-embedding $L \inj E$ corresponding to $\eta$. 
This is also a $K$-embedding. 
Therefore, $\Pm$ is true for $L/K$ and $m$. 
\end{proof}
To define filtrations of $G_K$, 
we show that the property $\Pm$ is stable under composition 
of finite Galois extensions of $K$ as follows:
\begin{proposition}\label{composite}
Let $L$ and $K'$ be finite Galois extensions of $K$. 
For a real number $m$, 
if $\Pm$ is true for both $L/K$ and $K'/K$, 
then $\Pm$ is also true for the composite extension $LK'/K$. 
In particular, we have $m_{LK'/K} \leq \max \{m_{L/K},m_{K'/K}\}$. 
\end{proposition}
\begin{proof}
Put $L':=LK'$. 
Assume $\Pm$ is true for $L/K$ and $K'/K$. 
Suppose there exists an $\cO_K$-algebra homomorphism 
$\eta:\cO_{L'} \to \cO_E/\fa_{E/K}^m$ for an algebraic extension $E$ of $K$. 
Then the composite maps defined by
\[\eta':\cO_L \inj \cO_{L'} \onto{\eta} \cO_E/\fa_{E/K}^m,\quad
\eta'':\cO_{K'} \inj \cO_{L'} \onto{\eta} \cO_E/\fa_{E/K}^m\]
are also $\cO_K$-algebra homomorphisms. 
Since $\Pm$ is true for both $L/K$ and $K'/K$, there exist 
$K$-embeddings $L \inj E$ and $K' \inj E$ corresponding to $\eta'$ and $\eta''$ 
respectively. 
Since $L/K$ and $K'/K$ are Galois extensions, 
we obtain a $K$-embedding $L' \inj E$. 
Therefore, $\Pm$ is true for $L'/K$. 
\end{proof}
By this proposition, we can define two increasing filtrations 
$(\bar{K}_{<m})_{m \in \mathbb{R}}$ and 
$(\bar{K}_{\leqslant m})_{m \in \mathbb{R}}$ 
of $\bar{K}$ as follows: For any real number $m$, 
$\bar{K}_{<m}$ (resp.\ $\bar{K}_{\leqslant m}$) 
is defined by the composite field 
of all finite Galois extensions $L$ of $K$ in $\bar{K}$ 
such that $m_{L/K} < m$ (resp.\ $m_{L/K} \leq m$). 
Then we put 
\[G_K^{\geqslant m}:=\Gal(\bar{K}/\bar{K}_{<m}),\quad 
G_K^{> m}:=\Gal(\bar{K}/\bar{K}_{\leqslant m}),\]
which are closed normal 
subgroups of $G_K$. 
Clearly, 
these subgroups 
form decreasing filtrations of $G_K$. 
\begin{remark}\label{RemImp}
In fact, Proposition \ref{pro1}, \ref{pro2} and \ref{composite} remain true 
in the case where the residue field $k$ may be imperfect, 
though we have to show 
the finiteness of  $m_{L/K}$ in Proposition \ref{pro1} by 
a different way via Proposition \ref{imperfect}. 
Hence the filtrations $G_K^{\geqslant m}$ and $G_K^{> m}$ 
can be defined even when the residue 
field of $K$ is imperfect. 
\end{remark}

The property $\Pm$ has the following property for 
unramified extensions of $K$: 
\begin{proposition}\label{etale}
Let $L$ be a finite Galois extension of $K$. 
Then the following conditions are equivalent:

\sn
$\mathrm{(i)}$ $L/K$ is unramified. 

\sn
$\mathrm{(ii)}$ $m_{L/K} \leq 0$. 

\sn
$\mathrm{(iii)}$ $m_{L/K}<1$. 
\end{proposition}
\begin{proof}
First, assume $L/K$ is unramified. 
Then we show that $\Pm$ is true for $L/K$ and $m>0$. 
Suppose there exists an $\cO_K$-algebra homomorphism 
$\eta:$ $\cO_L \to \cO_E/\fa_{E/K}^m$ for an algebraic extension 
$E$ of $K$. 
Since $L/K$ is unramified, 
$\eta$ lifts to 
an $\cO_K$-algebra homomorphism $\cO_L \to \cO_E$. 
Thus (i) implies (ii). 
Since it is clear that (ii) implies (iii), 
it is enough to verify that (iii) implies (i). 
To prove this, 
we show that if $L/K$ is not unramified, 
then $m_{L/K} \geq 1$. 
Let $K'$ be the maximal unramified subextension of $L/K$ and 
$\pi_K$ (resp.\ $\pi_L$) a uniformizer of $\cO_K$ (resp.\ $\cO_L$). 
Then there is an $\cO_K$-algebra homomorphism 
$\cO_L \to \cO_L/\pi_L \cO_L \cong \cO_{K'}/\pi_K \cO_{K'}=\cO_{K'}/\fa_{K'/K}^1$. 
However, there is no $K$-embedding $L \inj K'$. 
Hence (P$_1$) is not true for $L/K$, so that $m_{L/K} \geq 1$.  
\end{proof}
By the properties of the number $m_{L/K}$, 
our filtration $(G_K^{\geqslant m})_{m \in \mathbb{R}}$ has the 
following properties:
\begin{theorem}\label{properties} 
$\mathrm{(i)}$ For a real number $m \leq 0$, we have $G_K^{\geqslant m}=G_K$. 
Moreover, we have $\bigcap_m G_K^{\geqslant m}=1$ and 
$\overline{\bigcup_m G_K^{\geqslant m}}=G_K$. 

\noindent
$\mathrm{(ii)}$ Let $K'$ be a finite separable extension of $K$, 
of ramification index $e'$. 
We identify the Galois group $G_{K'}:=\Gal(\bar{K}/K')$ with a subgroup of $G_K$. 
Then,\ for a real number $m>0$, we have 
$G_{K'}^{\geqslant e'm} \subset G_K^{\geqslant m}$, 
with equality if $K'/K$ is unramified. 

\noindent
$\mathrm{(iii)}$ For a real number $0 < m \leq 1$, 
$G_K^{\geqslant m}$ is the inertia subgroup of $G_K$. 
\end{theorem}
\begin{proof}
The assertion (i) follows from Proposition $\ref{pro1}$. 
(iii) follows from Proposition $\ref{etale}$. 
The first assertion of (ii) follows from Proposition $\ref{pro2}$. 
Hence we prove the second assertion of (ii). 
Assume $K'/K$ is unramified. 
It suffices to show 
$\bar{K'}_{<m} \subset \bar{K}_{<m}$. 
By the definition of $\bar{K'}_{<m}$, 
it is enough to show that if a finite Galois extension $L'$ of $K'$ contained in $\bar{K}$ 
satisfies $m_{L^{'}/K^{'}}<m$, then $L' \subset \bar{K}_{<m}$. 
Since the case $L'=K'$ is true by Proposition \ref{etale}, 
we may assume $L' \not= K'$. 
Take the Galois closure $K''$ of $K'$ over $K$ in $\bar{K}$ 
and put $L'':= L'K''$. 
Note that $K''/K$ is 
an unramified Galois extension and $L''/K''$ is a Galois extension. 
Then we have
$m_{L^{''}/K^{''}} \leq m_{L^{'}/K^{'}} < m$ by Proposition \ref{pro2}. 
Let $\widetilde{L''}$ be the Galois closure of $L''$ 
over $K$ in $\bar{K}$. 
If $\widetilde{L''} = K''$, then Proposition \ref{etale} shows 
$m_{\widetilde{L''}/K} \leq 0 < m$, so that
$L' \subset \widetilde{L''} \subset \bar{K}_{<m}$. 
Thus we may assume $\widetilde{L''} \not= K''$. 
Any $\sigma \in \Gal(\widetilde{L''}/K)$ satisfies 
$\sigma(K'') = K''$ since $K''/K$ is a Galois extension, so that
\[m_{\sigma(L'')/K''} = m_{\sigma(L'')/\sigma(K'')} = m_{L^{''}/K^{''}} < m. \]
By this inequality and Proposition \ref{composite}, we have 
\[m_{\widetilde{L''}/K''} \leq \max \{ m_{\sigma(L'')/K''}\ |\ 
\sigma \in \Gal(\widetilde{L''}/K) \} < m \]  
since $\widetilde{L''}/K$ is the composite field of all the conjugate 
fields $\sigma(L'')$ $(\sigma \in \Gal(\widetilde{L''}/K))$.
Thus we have $m_{\widetilde{L''}/K} = m_{\widetilde{L''}/K''} < m$ 
by Proposition \ref{pro2} 
since $K''/K$ is unramified.
Therefore, we have $L' \subset \widetilde{L''} \subset \bar{K}_{<m}$. 
\end{proof}
\begin{remark}\label{WildInertia}
We can prove that $G_K^{>1}$ is the wild inertia subgroup of $G_K$ 
by using the property $\Pm$ together with the classical theory of 
Herbrand functions in a similar way to the proof of Proposition 
1.5, (ii), of \cite{Fon85}. 
However, we restricted ourselves here to showing what can be derived 
rather directly from $\Pm$. 
\end{remark}

\section{Ramification breaks}\label{RamificationBreaks}
\quad In this section, 
we compare our ramification filtration with the classical one. 
First, we recall the classical ramification theory for Galois 
extensions of $K$. 
Let $L$ be a finite Galois extension of $K$ with Galois group $G$. 
The order function $\ii_{L/K}$ is defined on $G$ by 
\[\ii_{L/K}(\sigma):=\inf_{a \in \cO_L}
 v_K(\sigma(a)-a),\ \sigma \in G.\]
Then the $i$th \emph{lower numbering ramification group} $G_{(i)}$ 
are defined for any real number $i$ by 
\[G_{(i)}:=\{ \sigma \in G\ |\ \ii_{L/K}(\sigma) \geq i \}.\]
The transition function $\ff_{L/K}:\mathbb{R} \to 
\mathbb{R}$ of $L/K$ is defined by 
\[\ff_{L/K}(u):=\int_0^u \sharp G_{(t)} dt\]
where $\sharp G_{(t)}$ is the cardinality of $G_{(t)}$. 
Then $\ff_{L/K}:\mathbb{R} \to \mathbb{R}$ is piecewise linear, strictly increasing and 
bijective (\cite{Serre}, Chap.\ IV, Sect.\ 3, Prop.\ 12). 
Denote by $\fp_{L/K}$ its inverse function. 
We also define another function $\uu_{L/K}$ on $G$ by 
\[\uu_{L/K}(\sigma):=\ff_{L/K}(\ii_{L/K}(\sigma)),\ \sigma \in G.\ \]
Then the $u$th \emph{upper numbering ramification group} $G^{(u)}$ 
are defined for any real number $u$ by 
\[G^{(u)}:=\{\sigma \in G\ |\ \uu_{L/K}(\sigma) \geq u \}.\]
For any non-negative real number $u$, 
we have $G^{(u)}=G^{u-1}$, where the latter is the $u$th upper numbering 
ramification group defined in \cite{Serre} (\emph{cf}.\ \cite{Fon85}, Rem.\ 1.2). 
We denote by $u_{L/K}$ (resp.\ $i_{L/K}$) 
the greatest upper (resp.\ lower) ramification break of $L/K$ 
defined by
\[u_{L/K}:=\inf \{u \in \mathbb{R}\ |\ G^{(u)}=1 \},\quad 
i_{L/K}:=\inf \{u \in \mathbb{R}\ |\ G_{(i)}=1 \}.\]
We put $u_{K/K}=-\infty$ by convention. 
The next lemma is a basic property of the number $u_{L/K}$:
\begin{lemma}\label{u_L/K}
For finite Galois extensions $M \subset L$ of $K$, 
we have $u_{M/K} \leq u_{L/K}$. 
\end{lemma}
\begin{proof}
By the compatibility with the quotient (\emph{cf}.\ \cite{Serre}, Chap.\ IV, 
Sect.\ 3, Prop.\ 14), 
$\Gal(L/K)^{(u)}=1$ implies $\Gal(M/K)^{(u)}=1$ for any real number $u$. 
Thus we obtain the inequality. 
\end{proof}
Fontaine proved the following proposition:
\begin{proposition}[\cite{Fon85}, Prop.\ 1.5]\label{Fontaine}
Let $L$ be a finite Galois extension of $K$ and $m$ a real 
number. 
Then there are the following relations:

\noindent
$\mathrm{(i)}$ If we have $m>u_{L/K}$, then $\Pm$ is true.  

\noindent
$\mathrm{(ii)}$ If $\Pm$ is true, then we have $m>u_{L/K}-e_{L/K}^{-1}$. 
\end{proposition}
By this proposition, 
we have the inequalities
\[u_{L/K}-e_{L/K}^{-1} \leq m_{L/K} \leq u_{L/K},\]
for a finite Galois extension $L$ of $K$. 
More precisely, 
we have the following equality:
\begin{proposition}\label{m=u}
For a finite Galois extension $L$ of $K$, 
we have $m_{L/K}=u_{L/K}$. 
\end{proposition}
\begin{proof}
It is enough to show that $\Pm$ is not true for $L/K$ and $m<u_{L/K}$. 
Suppose $L/K$ is unramified. 
Then we have $u_{L/K}=m_{L/K}=0$ if $L \not= K$, 
and $u_{L/K}=m_{L/K}=-\infty$ if $L=K$, so that the proposition follows. 
Therefore, we may assume $L/K$ is not unramified. 
The number $u_{L/K}$ is stable under unramified base change. 
Thus we may assume $L/K$ is a totally ramified extension by Proposition \ref{pro2}. 
If $L/K$ is a tamely ramified extension, 
$\Pm$ is not true even for $m=u_{L/K}=1$ because we can find a 
counter-example to $\Pm$ for $m=u_{L/K}$ as follows: 
Let $\pi_L$ (resp.\ $\pi_K$) be a unifirmizer of $\cO_L$ (resp.\ $\cO_K$). 
Then there is an $\cO_K$-algebra homomorphism 
$\cO_L \to \cO_L/\pi_L \cO_L \cong \cO_K/\pi_K \cO_K = \cO_K/\fa_{K/K}^1$. 
However, there is no $K$-embedding $L \inj K$. 
Therefore, we may assume $L/K$ is a wildly ramified 
extension. 
To prove this proposition, 
we shall find a counter-example to $\Pm$ for $L/K$ and 
$m=u_{L/K}-e'^{-1}$, where $e'$ can be taken to be an arbitrarily large number. 
Take a finite tamely ramified Galois extension $K'$ of $K$. 
Put $L':=LK'$ and $e':=e_{L'/K}$. 
If we apply (ii) of Proposition $\ref{Fontaine}$ to $L'/K$, 
then there exists an algebraic extension $E$ of $K$ such that 
there exists an $\cO_K$-algebra homomorphism $\eta:\cO_{L'} \to \cO_{E}/\fa_{E/K}^{m_0}$, 
but there is no $K$-embedding $L' \inj E$, 
where $m_0:=u_{L'/K} - e'^{-1}$. 
By Lemma $\ref{u_L/K}$, we have $m_0 \geq m_1$, 
where $m_1:=u_{L/K}-e'^{-1}$. 
Consider the two $\cO_K$-algebra homomorphisms defined by composite maps:
\[\eta':\cO_L \inj \cO_{L'} \onto{\eta} \cO_{E}/\fa_{E/K}^{m_0} 
\surj \cO_E/\fa_{E/K}^{m_1},\quad \eta'':\cO_{K'} \inj \cO_{L'} 
\onto{\eta} \cO_{E}/\fa_{E/K}^{m_0}.\]
Since $K'/K$ is a tamely ramified extension, 
we have $u_{K'/K} \leq 1$. 
On the other hand,\ since $L'/K$ is a wildly ramified extension, 
we have $e' m_0 > e'$ as shown in the proof of \cite{Fon85}, 
Proposition 1.5, (ii), hence we deduce $m_0 > 1$. 
Thus we have $m_0 > u_{K'/K}$. 
According to (i) of Proposition $\ref{Fontaine}$ for $K'/K$,  
there exists a $K$-embedding $K' \inj E$ corresponding to $\eta''$. 
If we suppose there exists a $K$-embedding $L \inj E$, 
then there exists a $K$-embedding $L'=LK' \inj E$ since $L/K$ and 
$K'/K$ are Galois extensions. 
This is a contradiction. 
Therefore, $\Pm$ is not true for $L/K$ and $m=m_1$.
Hence the result follows. 
\end{proof}
\begin{remark}
By Proposition $\ref{composite}$, Lemma $\ref{u_L/K}$ and Proposition 
$\ref{m=u}$, 
we deduce the equality $u_{LK'/K}=\max \{u_{L/K},u_{K'/K}\}$ for any 
finite Galois extensions $L$ and $K'$ of $K$. 
\end{remark}
\begin{remark}\label{monogenic}
In the above proposition, 
we proved the equality $m_{L/K}=u_{L/K}$ with 
the assumption that the residue filed $k$ is perfect. 
We are also interested in the case where $k$ may be imperfect. 
In Chapter IV of \cite{Serre}, the ramification filtration is defined 
in the case where $L/K$ is unferociously\footnote{
We mean by an \emph{unferociously} ramified extension 
$L/K$ a finite algebraic extension whose residue field extension 
is separable.} 
ramified. 
Our proof of Proposition \ref{m=u} remains true in this case
since so does Fontaine's proof of Proposition \ref{Fontaine} and
the composite field of $L/K$  and any tamely ramified 
extension $K'/K$  is still unferociously ramified.
\end{remark}
Theorem $\ref{Main}$ follows from Propositions $\ref{m=u}$ and $\ref{filtration}$.

\section{The ramification theory of Abbes and Saito}\label{ImperfectSection}
\quad First, we recall the ramification theory of Abbes and Saito 
(\cite{Abbes-Saito 1}, \cite{Abbes-Saito 2}).
In Subsection $\ref{ImperfectPm}$,
we generalize the property $\Pm$ to the imperfect residue field case.
In Subsection $\ref{Translate}$,
we translate our results in Section $\ref{RamificationBreaks}$ into the language 
of the ramification theory of Abbes and Saito.
Let $K$ be a complete discrete valuation field 
whose residue field $k$ may not be perfect. 
Let $K^{\mathrm{alg}}$ be a fixed algebraic closure of $K$, 
$\bar{K}$ the separable closure of $K$ in $K^{\mathrm{alg}}$ and $G_K:=\Gal(\bar{K}/K)$ 
the absolute Galois group. 
Abbes and Saito defined a decreasing filtration 
$(G_K^m)_{m \geq 0}$ by closed normal subgroups 
$G_K^m$ of $G_K$ indexed by 
rational numbers $m \geq 0$, in such a way 
that $\bigcap_{m \geq 0} G_K^m=1,\ G_K^0=G_K$ 
and $G_K^1$ is the inertia subgroup of $G_K$. 
It is defined by using certain functors $F$ and $F^m$ from the category 
$\fe_K$ of finite \'etale $K$-algebras to the category $\bs_K$ 
of finite $G_K$-sets. 
We recall the definition of $F$ and its quotients $F^m$ 
for positive rational numbers $m$. 
Let $L$ be a finite \'etale $K$-algebra and $\cO_L$ the integral 
closure of $\cO_K$ in $L$. 
We define 
$F(L):=\mathrm{Hom}_K(L,\bar{K})=
\mathrm{Hom}_{\cO_K}(\cO_L,\cO_{\bar{K}})$. 
The functor $F$ gives an 
anti-equivalence 
of $\fe_K$ with $\bs_K$, thereby making 
$\fe_K$ a Galois category. To define $F^m$, we proceed as follows: 
An \emph{embedding} of $\cO_L$ is a pair $(\mathbb{B}, \mathbb{B} \to 
\cO_L)$ consisting of an $\cO_K$-algebra $\mathbb{B}$ 
which is formally of finite type\footnote{
We say that an $\cO_K$-algebra $A$ is \emph{formally of finite type} 
over $\cO_K$ if $A$ is semi-local, $\mathfrak{m}_A$-adically complete 
Notherian and the quotient $A/\mathfrak{m}_A$ is finite over $k$, 
where $\mathfrak{m}_A$ is the radical of $A$ 
(\emph{cf}.\ \cite{Abbes-Saito 2}, Sect.\ 1). } 
and formally smooth over $\cO_K$ and a surjection $\mathbb{B} \to \cO_L$ of $\cO_K$-algebras 
which induces an isomorphism 
$\mathbb{B}/\mathfrak{m}_{\mathbb{B}} \to \cO_L/\mathfrak{m}_L$, where $\mathfrak{m}_{\mathbb{B}}$ and 
$\mathfrak{m}_L$ are respectively the radicals of $\mathbb{B}$ and $\cO_L$ 
($cf$.\ \cite{Abbes-Saito 2}, Def.\ 1.1). 
Let $I$ be the kernel of the surjection $\mathbb{B} \to \cO_L$. 
Write $m = m_2/m_1$ for some positive integers $m_1$ and $m_2$. 
Then the affinoid algebra $\mathbb{B}
[I^{m_1}/\pi_K^{m_2}]^{\land} \otimes_{\cO_K} K$ 
does not depend on 
the presentation of $m$ (\cite{Abbes-Saito 2}, Lem.\ 1.4, 4), 
where $\pi_K$ is a 
uniformizer of $K$ and $\land$ means the $\pi_K$-adic completion. 
Hence we denote this ring by $\ssb^m$. 
Let $X^m(\mathbb{B} \to \cO_L)$ be the affinoid 
variety $\mathrm{Sp}(\ssb^m)$ associated with $\ssb^m$. 
For any affinoid variety $X$ over $K$, let $\pi_0(X_{\bar{K}})$ 
denote the set $\underset{\longleftarrow}{\lim}_{K'} \pi_0(X \otimes_K K')$ 
of geometric connected components, where $K'$ runs through the finite separable 
extensions of $K$ contained in $\bar{K}$. 
Then we define the functor $F^m$ by 
\[F^m(L):=\lim_{\underset{(\mathbb{B} \to \cO_L)}{\longleftarrow}} 
\pi_0(X^m(\mathbb{B} \to \cO_L)_{\bar{K}}),\]
where $(\mathbb{B} \to \cO_L)$ runs through the category of embeddings of $\cO_L$ 
(\emph{cf}.\ \cite{Abbes-Saito 2}, Def.\ 1.1).
The projective system in the right-hand 
side is constant (\cite{Abbes-Saito 2}, Lem.\ 1.9).
The finite set $F(L)$ can be identified with a subset of 
$X^m(\mathbb{B} \to \cO_L)(\bar{K})$, and this induces 
a natural surjective map $F(L) \to F^m(L)$. 
The $m$th ramification subgroup $G_K^m$ is characterized by the property that 
$F(L)/G_K^m=F^m(L)$ for all $L$. 
If the residue field of $K$ is perfect, 
this filtration $(G_K^m)_m$ defined as above coincides with 
the classical one $(G_K^{(m)})_m$ defined in Section 3 
(\emph{cf}.\ \cite{Abbes-Saito 1}, Subsect.\ 6.1).

\subsection{Generalization of Fontaine's proposition}\label{ImperfectPm}
\quad In this subsection, 
we generalize Fontaine's proposition to the imperfect residue field case. 
Let $L$ be a finite Galois extension of $K$ and $m$ a positive rational 
number. 
We define the property $\Pm$ and the number 
$m_{L/K}$ in the same way 
as those in the Introduction. 
For an affinoid variety $X$ over $K$ and a point $x \in X(K^{\mathrm{alg}})$, 
we denote by $K(x)$ and $X_x$, respectively, the definition field of $x$, and 
the geometric connected component of $X$ which contains $x$. 
The ring $\cO_L$ is a complete intersection over $\cO_K$. 
Namely, we have $\cO_L \cong \cO_K[T_1,\dots,T_n]/(f_1,\dots,f_n)$ 
(\cite{Abbes-Saito 1}, Lem.\ 7.1). 
We denote by $z_1,\dots,z_d$ the common zeros of $f_1,\dots,f_n$ 
in $\bar{K}^n$. 
Let $I:=(f_1,\dots,f_n)$ be the ideal of $\cO_K[T_1,\dots,T_n]$ generated by 
$f_1,\dots,f_n$. 
Consider the surjection $\varphi: \cO_K[T_1,\dots,T_n] \to \cO_K[T_1,\dots,T_n]/(f_1,\dots,f_n) 
\cong \cO_L$. 
Put $w_i:=\varphi(T_i)$ for $i=1,\dots,n$. 
Then the formal completion $\mathbb{B} \to \cO_L$ of 
$\cO_K[T_1,\dots,T_n] \to \cO_L$, 
where 
$\mathbb{B}:=\underset{\longleftarrow r}{\lim} \cO_K[T_1,\dots,T_n]/
I^r$, is an embedding of $\cO_L$. 
Let $X^m:=X^m(\mathbb{B} \to \cO_L)$ be the affinoid variety over $K$ 
associated with this embedding. 
Then we have 
\[X^m(K^{\mathrm{alg}})=\{ x \in \cO_{K^{\mathrm{alg}}}^n\ |\ 
v_K(f_i(x)) \geq m\ (i=1,\dots,n)\}.\]
\begin{remark}
If $m$ is not a rational number, then $X^m$ does not form a $K$-affinoid 
variety. 
\end{remark}
\begin{lemma}\label{hom}
Let $m$ be a positive rational number and $E/K$ an algebraic extension. 
Then, the map $X^m(E) \to \mathrm{Hom}_{\cO_K}
(\cO_L,\cO_E/\fa_{E/K}^m)$ sending $(x_1,\dots,x_n)$ to the homomorphism 
defined by $w_i \mapsto x_i \pmod{\fa_{E/K}^m}$, 
is surjective. 
\end{lemma}
\begin{proof}
Obvious. 
\end{proof}
Consider the following property for $L/K$ and $m$:
\begin{quote}
\begin{itemize}
\item[$\PPm$] \emph{For any $x \in X^m(K^{\mathrm{alg}})$, 
there exists a common zero $z$ of $f_1,\dots,f_n$ in $\bar{K}^n$ 
which is $K(x)$-rational}. 
\end{itemize}
\end{quote}
We can easily check that $\PPm$ is equivalent to $\Pm$ if $m$ is a 
positive rational number by Lemma \ref{hom}. 
On the other hand, we consider the following property for $L/K$ and $m$:
\begin{quote}
\begin{itemize}
\item[$\QQm$] \emph{For any $x \in X^m(K^{\mathrm{alg}})$, 
there exists a common zero $z$ of $f_1,\dots,f_n$ in $\bar{K}^n$ 
such that $x \not\in X^m_{z_i}$ for any $z_i$ except $z$}.
\end{itemize}
\end{quote}
By definition, the property $\QQm$ is equivalent to the bijectivity of $F(L) \to F^m(L)$.\
Let $c_{L/K}$ be the conductor of $L/K$ (\cite{Abbes-Saito 1}, Def.\ 6.3), 
which is defined by 
\begin{center}
$c_{L/K}:= 
\inf \{m \in \mathbb{Q}_{\geq 0}\ |\ \QQm$ is true for $L/K.\}$.
\end{center}
If the residue field $k$ is perfect, 
we have $c_{L/K}=u_{L/K}$ for any finite Galois extension $L$ of $K$. 
We can show the following proposition which is a generalization of 
(i) of Proposition $\ref{Fontaine}$ to the imperfect residue field case:
\begin{proposition}\label{imperfect}
If $\QQm$ is true, then $\PPm$ is true. 
In particular, we have the inequality 
$m_{L/K} \leq c_{L/K}$. 
\end{proposition}
\begin{proof}
We need the following lemma 
which is a version of Krasner's lemma. 
This is due to Hiranouchi and Taguchi. 
\begin{lemma}\label{KrasnerImperfect}
Let $X$ be an affinoid variety over $K$. 
Let $x \in X(\bar{K})$ and $y \in X(K^{\mathrm{alg}})$. 
Assume that any $G_K$-conjugates of $x$ different from $x$ is not contained in 
$X_x$ and that $y$ is in the geometric connected component 
$X_x$. 
Then $K(x) \subset K(y)$. 
\end{lemma}
This lemma is proved in the same way as the classical one. 
\begin{proof}
If $\sigma$ $\in$ $\mathrm{Hom}_{K(y)}(K(x,y),K^{\mathrm{alg}})$, 
we have $y \in X_{\sigma(x)}$ and $y \in X_{x}$, 
so that we have $X_{\sigma(x)}=X_x$. 
Hence $\sigma$ fixes $x$ by the assumption on $x$. 
Thus we have $K(x)$ $\subset$ $K(y)$. 
\end{proof}
Now we can prove Proposition \ref{imperfect} as follows: 
Let $m$ be a positive rational number and $x$ a point of 
$X^m(K^{\mathrm{alg}})$. 
By the property $\QQm$, 
there exists a zero $z$ of $f_1,\dots,f_n$ such that $x \in X^m_z$ but 
$x \not\in X^m_{z_i}$ for any $z_i \ne z$. 
Then we have $z_i \not\in X_z^m$ for any $z_i \not= z$. 
Indeed, if $z_i \in X_z^m$, then $X_{z_i}^m = X_z^m$, 
which contradicts $x \in X^m_z$ and $x \not\in X_{z_i}^m$ for any $z_i \not= z$. 
Thus $K(z) \subset K(x)$ by Lemma \ref{KrasnerImperfect}. 
Hence $\PPm$ is true. 
\end{proof}
\begin{remark}
The author does not know whether the equality $m_{L/K}=c_{L/K}$ remains true 
in the case where the residue field of $K$ is imperfect. 
However, we can show at least the following: 
\end{remark}
\begin{proposition}
Let $L$ be a finite Galois extension of $K$. 
Let $K'$ be a weakly unramified\footnote{
We mean by a \emph{weakly unramified} extension 
$K'/K$ a finite algebraic extension such that $e_{K'/K}=1$.} extension of $K$ such that 
$L':=LK'/K'$ is unferociously ramified 
(the existence of such an extension is proved in 
\cite{Abbes-Saito 1}, Append.\ Cor.\ A.2). 
Then we have $c_{L'/K'} \leq m_{L/K}$. 
\end{proposition}
\begin{proof}
We have $m_{L'/K'} \leq m_{L/K}$ by Proposition \ref{pro2} 
(\emph{cf}.\ Rem.\ \ref{RemImp}). 
Since $L'/K'$ is unferociously ramified, 
we can apply Proposition \ref{m=u} to $L'/K'$ (\emph{cf}.\ Rem.\ \ref{monogenic}). 
Then we have $c_{L'/K'} = m_{L'/K'}$. 
Thus the desired inequality $c_{L'/K'} \leq m_{L/K}$ holds. 
\end{proof}

\subsection{Comparison with the ramification theory of Abbes and Saito}\label{Translate}
\quad In this subsection, 
we translate our results in Section $\ref{RamificationBreaks}$ 
into the language of the ramification theory of Abbes and Saito. 
Let $K$ be a complete discrete valuation field with perfect residue field and $L$ a 
finite Galois extension of $K$.  
We define a non-Archimedean valuation on $K^{\mathrm{alg}}$ by $|y|=\theta^{v_K(y)}$, 
where $0<\theta<1$ is a real number. 
Fix a generator $z$ of $\cO_L$ as an $\cO_K$-algebra.  
Let $P$ be the minimal polynomial of $z$ over $K$, and 
$z=z_1,\dots,z_d$ the zeros of $P$ in $\bar{K}$. 
Let $X^m$ be the affinoid variety over $K$ as defined in the 
previous subsection, 
so that $X^m(K^{\mathrm{alg}})$ $=$ 
$\{x \in \cO_{K^{\mathrm{alg}}}\ |\ v_K(P(x)) \geq m \}$. 
If the residue field of $K$ is perfect, 
we can rewrite $\PPm$ for $L/K$ and a positive 
rational number $m$ as follows: 
\begin{quote}
\begin{itemize}
\item[$\PPPm$] \emph{For any $x \in X^m(K^{\mathrm{alg}})$,  
there exists a zero $z$ of $P$ in $\bar{K}$ 
which is $K(x)$-rational}. 
\end{itemize}
\end{quote}
On the other hand, we consider the following property: 
\begin{quote}
\begin{itemize}
\item[$\Qm$] \emph{For any $x \in X^m(K^{\mathrm{alg}})$, 
there exists a zero $z$ of $P$ in $\bar{K}$ such that 
$|z - x| = \min_i |x-z_i|$ and $|z - x| <\min_{i \not= 1 } |z-z_i|$}. 
\end{itemize}
\end{quote}
\begin{proposition}\label{QmEquiv}
The properties $\QQm$ and $\Qm$ are equivalent. 
\end{proposition}
\begin{proof}
Put $D(z_i,\theta^m):=\{ x \in \cO_{K^{\mathrm{alg}}}\ |\ |x-z_i| \leq \theta^m \}$ for 
$i=1,\dots,d$. 
The disc $D(z_i,\theta^m)$ is connected and contains $z_i$. 
Denote $\fp:=\fp_{L/K}$ for simplicity. 
Then we have the following by Lemma \ref{Herbrand} below: 
\[
\begin{split}
X^m(K^{\mathrm{alg}}) &= \{ x \in \cO_{K^{\mathrm{alg}}}\ |\ |P(x)| \leq \theta^m \} \\
&= \{ x \in \cO_{K^{\mathrm{alg}}}\ |\ \min_i |x-z_i| \leq \theta^{\fp(m)} \} \\
&= \underset{i}{\bigcup}\ D(z_i,\theta^{\fp(m)}). 
\end{split}
\]
The property $\QQm$ is true if and only if 
$X_{z_i}^m$ $(1 \leq i \leq d)$ 
are disjoint. 
The assertion of the proposition 
follows from the following equivalences: 
\[
\begin{split}
X_{z_i}^m \cap X_{z_j}^m = \emptyset \quad (i \not= j) &\iff 
\underset{i \not= j}{\min}\ |z_i - z_j| > \theta^{\fp(m)} \\ 
&\iff \min_i |z_i - x| < \min_{i \not= j} |z_i - z_j|\ 
\mathrm{for}\ \mathrm{all}\ x \in X^m(K^{\mathrm{alg}}).
\end{split}
\]
The first equivalence is proved as follows: 
Let $i \not= j$. 
Assume $X_{z_i}^m$ $(1 \leq i \leq d)$ are disjoint.  
Then we have $z_j \not\in X_{z_i}^m$. 
On the other hand, we have $D(z_i,\theta^{\fp(m)}) \subset X_{z_i}^m$ 
since $D(z_i,\theta^{\fp(m)})$ is connected and contained in $X^m(K^{\mathrm{alg}})$. 
Hence the zero $z_j$ is not contained in $D(z_i,\theta^{\fp(m)})$, so that 
$|z_i -z_j| >\theta^{\fp(m)}$. 
Conversely, suppose $|z_i -z_j| >\theta^{\fp(m)}$. 
Then we have $D(z_i,\theta^{\fp(m)}) \cap D(z_j,\theta^{\fp(m)}) = \emptyset$. 
Hence we obtain the decomposition 
$X^m(K^{\mathrm{alg}}) = \bigsqcup_i D(z_i,\theta^{\fp(m)})$. 
Thus we deduce $X_{z_i}^m = D(z_i,\theta^{\fp(m)})$. 
In particular, the connected components $X_{z_i}^m$ $(0 \leq i \leq d)$ are disjoint. 
Finally, we prove the second equivalence. 
Assume $\underset{i \not= j}{\min} |z_i - z_j| > \theta^{\fp(m)}$. 
For a point $x \in X^m(K^{\mathrm{alg}})$, 
take $i_1$ such that $x \in D(z_{i_1},\theta^{\fp(m)})$. 
Then we have  
\[ \min_i |z_i - x| \leq |z_{i_1} - x| \leq \theta^{\fp(m)} < \min_{i \not= j} |z_i - z_j|. \]
Conversely, assume $\min_i |z_i - x| < \min_{i \not= j} |z_i - z_j|$ for 
any $x \in X^m(K^{\mathrm{alg}})$.
Take $y \in K^{\mathrm{alg}}$ such that $|y| = \theta^m$ $(<1)$ 
and take $x \in K^{\mathrm{alg}}$ such that $P(x)=y$. 
Then we have $x \in \cO_{K^{\mathrm{alg}}}$ and $|P(x)| = \theta^m$. 
In particular, this shows $x \in X^m(K^{\mathrm{alg}})$. 
By assumption and Lemma \ref{Herbrand} below, we have the inequality
\[ \theta^{\fp(m)} = \min_i |z_i - x| < \min_{i \not= j} |z_i - z_j|. \]

\end{proof}
\begin{lemma}[\cite{Fon85}, Prop.\ 1.4]\label{Herbrand}
Let $x$ be an element of $K^{\mathrm{alg}}$. 
Put $i := \sup_i v_K(z_i - x)$ and 
$u:= v_K(P(x))$. 
Then we have
\[ u = \ff_{L/K}(i),\quad \fp_{L/K}(u)=i. \]
\end{lemma}
We obtain the following consequences:
\begin{proposition}
We have the following relations:

\noindent
$\mathrm{(i)}$ If $\Qm$ is true, then $\PPPm$ is true. 

\noindent
$\mathrm{(ii)}$ If $\PPPm$ is true, 
then $(\mathrm{Q}''_{m+\varepsilon})$ is true 
for any $\varepsilon>0$. 

\noindent
In particular, we have the equality 
$m_{L/K}=c_{L/K}$. 
\end{proposition}
\begin{proof}
The above (i) is the special case of Proposition $\ref{imperfect}$. 
(ii) follows from Propositions \ref{m=u}, Proposition \ref{QmEquiv} and the equality 
$u_{L/K}=c_{L/K}$. 
\end{proof}

\section{Appendix}
In this section, 
we prove a Galois theoretic property of a 
filtration of the absolute Galois group of an arbitrary field. 
This section is independent of the other sections. 
Let $K$ be a field, $\bar{K}$ a fixed separable closure of $K$, 
$G_K:=\Gal(\bar{K}/K)$ the absolute Galois group of $K$ and 
$\mathfrak{G}$ the set of all finite Galois extensions of $K$ contained in $\bar{K}$. 
Throughout this Appendix, 
all separable extensions of $K$ are assumed to be subfields of $\bar{K}$. 
Let $R$ be a totally ordered set. 

\begin{definition}\label{compatibility}
Assume we are given a system of decreasing filtrations 
$(\Gal(L/K)^u)_{L/K \in \mathfrak{G}, u \in R}$. 
Then we say that the system of filtrations $(\Gal(L/K)^u)_{L/K \in \mathfrak{G}, u \in R}$ is 
\emph{quotient-compatible} if, 
for any $L$, $L'$ $\in \mathfrak{G}$ such that 
$L' \subset L$, the image of $\Gal(L/K)^u$ under the 
natural projection $\Gal(L/K) \to \Gal(L'/K)$ coincides with 
$\Gal(L'/K)^u$. 
\end{definition}
\begin{proposition}
There is a natural one-to-one correspondence between 
the set of decreasing filtrations $(G_K^u)_{u \in R}$ on $G_K$ consisting of 
closed subgroups of $G_K$ and the set of quotient-compatible systems 
of decreasing filtrations $(\Gal(L/K)^u)_{L/K \in \mathfrak{G}, u \in R}$. 
\end{proposition}
\begin{proof}
Assume we are given a decreasing filtration $(G_K^u)_{u \in R}$ on $G_K$ 
consisting of closed subgroups of $G_K$. 
Let $L$ be a finite Galois extension of $K$ with Galois group $G$. 
Then a decreasing 
filtration $G^u$ can be defined by the image of $G_K^u$ by 
the restriction map $G_K \to G$. 
Conversely, suppose we are given a quotient-compatible system of decreasing filtrations 
$(\Gal(L/K)^u)_{L/K \in \mathfrak{G}, u \in R}$. 
For any finite Galois extensions $L' \subset L$ of $K$, 
the compatibility with the quotient induces a natural projection 
$\Gal(L/K)^u \to \Gal(L'/K)^u$ by the restriction map. 
Hence we can define a decreasing filtration 
$G_K^u$ on $G_K$ by 
\[G_K^u:=\underset{\longleftarrow}{\lim}\ \Gal(L/K)^u,\]
where $L$ runs through the set of all finite Galois extensions of $K$ 
contained in $\bar{K}$. 
This correspondence induces the desired bijection. 
\end{proof}
\begin{definition}
Let $G$ be a set and $(G^u)_{u \in R}$ a decreasing filtration on $G$. 
Then we say that $(G^u)_{u \in R}$ is \emph{separated} if 
$\bigcap_u G^u=1$ 
and $(G^u)_{u \in R}$ is \emph{left continuous} if 
$G^u=\bigcap_{m<u} G^m$. 
\end{definition}
Let $(G_K^u)_{u \in R}$ be a decreasing filtration on $G_K$ 
which is separated and left continuous, and 
$L$ a finite Galois extension of $K$ with Galois group $G$. 
Put $G_K^{u+}:=\overline{\bigcup_{u'>u} G_K^{u'}}$, 
where the overline means the closure with respect to Krull topology. 
Then we denote by $\bar{K}_{(u)}$ (resp.\ $\bar{K}_{(u+)}$) 
the fixed field of $\bar{K}$ by $G_K^u$ (resp.\ $G_K^{(u+)}$). 
Define $G^u$ (resp.\ $G^{u+}$) 
as the image of $G_K^u$ (resp.\ $G_K^{u+}$) by the restriction map 
$\pi:G_K \to G$. 
Put 
\[u_{L/K}:=\inf \{ u \in R\ |\ G^u=1 \},\]
assuming that the infimum exists in $R$. 
We denote by $\bar{K}_{<u}$ (resp.\ $\bar{K}_{\leqslant u}$) 
the union of all finite Galois extension 
$L$ of $K$ in $\bar{K}$ such that $u_{L/K}<u$ (resp.\ $u_{L/K} \leq u$). 
\begin{proposition}\label{filtration}
We have $\bar{K}_{<u}=\bar{K}_{(u)}$ and 
$\bar{K}_{\leqslant u}=\bar{K}_{(u+)}$ for any $u \in R$. 
\end{proposition}
\begin{proof}
If $L$ is a finite Galois extension of $K$ with Galois group $G$, 
then the left continuousness makes $G^{u_{L/K}} \ne 1$. 
Hence $u_{L/K}<u$ (resp.\ $u_{L/K} \leq u$) 
is equivalent to $G^u=1$ (resp.\ $G^{u+}=1$). 
This is equivalent to $G_K^u \subset \mathrm{Ker}(\pi)=G_L$ 
(resp.\ $G_K^{u+} \subset G_L$). 
The result follows it. 
\end{proof}


\end{document}